\documentclass[reqno,11pt]{amsart}

\newtheorem{rem}{Remark}

\newcommand\eps\varepsilon
\newcommand\ph\varphi
\newcommand\kap\varkappa

\usepackage{enumerate}
\usepackage{graphicx}
\usepackage{float}

\begin{document}

\title[On Periodic Solutions to Some Lagrangian System]
{On Periodic Solutions to Some Lagrangian System With Two Degrees of Freedom}

\author[Oleg Zubelevich]{Oleg Zubelevich\\ \\\tt
 Dept. of Theoretical mechanics,  \\
Mechanics and Mathematics Faculty,\\
M. V. Lomonosov Moscow State University\\
Russia, 119899, Moscow,  MGU \\
 }
\date{}
\thanks{Partially supported by grants
 RFBR  12-01-00441.,  Science Sch.-2964.2014.1}
\subjclass[2000]{34C25, 	70F20     }
\keywords{Lagrangian systems, periodic solutions.}

\begin{abstract}A Lagrangian system with two degrees of freedom is considered. The configuration space of the system is a cylinder. A large class of periodic solutions has been found. The solutions are not homotopy equivalent  to each other.
\end{abstract}

\maketitle
\numberwithin{equation}{section}
\newtheorem{theorem}{Theorem}[section]
\newtheorem{lemma}[theorem]{Lemma}
\newtheorem{definition}{Definition}[section]

\section{Statement of the Problem and Main Result}
This short note is devoted to  the following dynamical system.

\begin{figure}[H]
\centering
\includegraphics[width=45mm, height=30 mm]{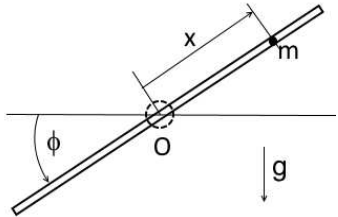}
\caption{the tube and the ball \label{overflow}}
\end{figure}

{\it A thin  tube can rotate freely in the vertical plane about a fixed horizontal axis passing through its centre $O$. A moment of inertia of the tube about this axis is equal to $J$. The mass of the tube is distributed symmetrically such that tube's centre of mass is placed at the point $O$.

Inside the tube there is a small ball  which can slide  without friction. The mass of the ball is $m$. The  ball can pass by the point $O$ and fall out from the ends of the tube.

The system  undergoes  the standard gravity field $\boldsymbol g$.}

It seems to be evident that for typical motion the ball reaches an end of the tube and falls down out the tube. It is surprisingly, at least for the first glance, that this system has very many periodic solutions such that the tube turns around several times during the period. 

The sense of generalised coordinates $\phi,x$ is clear from  Figure \ref{overflow}.

A kinetic energy and a potential of the system are given by the formulas
$$T=\frac{1}{2}\Big(mx^2+J\Big)\dot\phi^2+\frac{1}{2}m\dot x^2,\quad V=mgx\sin\phi.$$
By the suitable choice of dimension of units we obtain $$J=1,\quad g=1,\quad m=1.$$ So that a Lagrangian of the system is
\begin{equation}\label{Lagr}L(x,\phi,\dot x,\dot \phi)= \frac{1}{2}\Big(x^2+1\Big)\dot\phi^2+\frac{1}{2}\dot x^2 -x\sin\phi.\end{equation}

\begin{theorem}\label{dt5y555}
For any constants $\omega>0,\quad k\in\mathbb{N}$ system (\ref{Lagr}) has a solution $\phi(t),x(t),\quad t\in\mathbb{R}$ such that 

1) $x(t)=-x(-t),\quad  \phi(t)=-\phi(-t);$

2) $x(t+\omega)=x(t),\quad \phi(t+\omega)=\phi(t)+2\pi k.$
\end{theorem}
This theorem means that if $\omega$ and $k$ are given and the tube is long enough then the system has an $\omega-$periodic motion and the tube turns around  $k$ times during the period.

\section{ Proof of Theorem \ref{dt5y555}}
\subsection{Preliminary Remarks}
Introduce a space
$$H^1_o(-a,a)=\{u\in H^1(-a,a)\mid u(-t)=-u(t)\},\quad a\in(0,\infty).$$
  Recall that the Sobolev space  $H^1(-a,a)$ is compactly embedded to $C[-a,a]$.

\begin{lemma}\label{xddddd}
Let $u\in H^1_o(-a,a)$ then the following inequalities hold
$$\|u\|^2_{L^2(0,a)}\le \frac{a^2}{2}\|\dot u\|^2_{L^2(0,a)},\quad \|u\|^2_{C[0,a]}\le a \|\dot u\|^2_{L^2(0,a)}.$$\end{lemma}
This Lemma is absolutely  standard, we bring its proof just for completeness of exposition. \begin{rem}  Lemma \ref{xddddd} implies that the function $u\mapsto \|\dot u\|_{L^2(0,a)}$ is a norm  of $H^1_o(-a,a)$ and this norm is equivalent to the standard norm of $H^1(-a,a)$.\end{rem}

{\it Proof of Lemma \ref{xddddd}.}
We prove only the first inequality the second one goes in the same way. First assume that a function $u\in H^1(-a,a)$ is smooth.  From the formula
$$u(t)=\int_0^t\dot u(s)ds$$ it follows that
$$\int_0^au^2(s)ds=\int_0^a\Big(\int_0^t\dot u(s)ds\Big)^2dt.$$
It remains to observe that by the Cauchy inequality
$$\Big|\int_0^t\dot u(s)ds\Big|\le \int_0^t|\dot u(s)|ds\le \|\dot u\|_{L^2(0,a)}\Big(\int_0^tds\Big)^{1/2},\quad t\in[0,a].$$
Since the space of smooth functions is dense in $H^1(-a,a)$, the inequality under consideration holds for all $u\in H^1(-a,a)$.

The Lemma is proved.

\begin{lemma}\label{dfdfff}Being endowed with a collection of seminorms 
\begin{equation}\label{sdfrgcc}\|u\|_n=\|u\|_{H^1(-n,n)},\quad n\in\mathbb{N}\end{equation} the space $H^1_{\mathrm{loc}}(\mathbb{R})$ turns to  a reflexive Fr\'echet space.\end{lemma}
\begin{rem}\label{dfgdggg}It would be  more accurate  to write formula (\ref{sdfrgcc}) as follows
$\|u\|_n=\|u\mid_{[-n,n]}\|_{H^1(-n,n)}$, where $\mid_{[-n,n]}$ is the operation of restriction to the interval $[-n,n]$. Nevertheless here and in the sequel we will hold  this little bit informal notation. It will not generate a misleading.\end{rem}
Surely Lemma \ref{dfdfff} is  a trivial and well-known fact. Nevertheless, we did not encounter it in the textbooks, so we present its proof.

{\it Proof of Lemma \ref{dfdfff}.} It is clear that the space $H^1_{\mathrm{loc}}(\mathbb{R})$ is compete, thus it is a barrelled space 
\cite{iosida}. 

The space  $H^1_{\mathrm{loc}}(\mathbb{R})$ is a projective limit of the spaces $H^1(-n,n)$ with respect to the restriction operators $$H^1_{\mathrm{loc}}(\mathbb{R})\to H^1(-n,n).$$

The projective limit of  reflexive  spaces is a semi-reflexive  space 
\cite{Shefer}.
 A barreled semi-reflexive space is a reflexive space
 \cite{iosida}. 
Consequently, $H^1_{\mathrm{loc}}(\mathbb{R})$ is a reflexive space.

The Lemma is proved.

Determine the following subspaces
$$H^1_o(\mathbb{R})=\{u\in H^1_{\mathrm{loc}}(\mathbb{R})\mid u(-t)=-u(t)\}$$ and
 $$X_\omega=\{x\in H^1_o(\mathbb{R})\mid x(t+\omega)=x(t)\}.$$
They both are closed. Moreover, from Lemma \ref{xddddd} it follows that a function $x\mapsto \|\dot x\|_{L^2(0,\omega)}$ is a norm in $X_\omega$ and the topology of this norm is equivalent to the one inherited of  $H^1_{\mathrm{loc}}(\mathbb{R})$. So $X_\omega$ is a Banach space.

\begin{lemma}\label{ssdddd} The spaces $H^1_o(\mathbb{R}),X_\omega$ are  reflexive.\end{lemma}
The proof of this lemma almost literally repeats the proof of Lemma \ref{dfdfff}. Just note that the space
$H_o^1(-n,n)$ is a reflexive space because it is a real Hilbert space with standard inner product
$$(u,v)=\int_{(-n,n)}u(t)v(t)dt+\int_{(-n,n)}\dot u(t)\dot v(t)dt.$$

Introduce a set
$$\Phi_{k,\omega}=\{\phi\in H^1_o(\mathbb{R})\mid\phi(t+\omega)=\phi(t)+2\pi k\}.$$
The set $\Phi_{k,\omega}$ is closed and convex  in $H^1_o(\mathbb{R})$.

With the help of Lemma \ref{xddddd} it is not hard to show that the function $\rho(u,v)=\|\dot u-\dot v\|_{L^2(0,\omega)}$ determines a metric on $\Phi_{k,\omega}$ and this metric endows $\Phi_{k,\omega}$ with the same topology   as the space $H^1_o(\mathbb{R})$ does.

\subsection{The Action}Our goal is to prove that  a functional
$$S:X_\omega\times H^1_o(\mathbb{R})\to\mathbb{R},\quad S(x,\phi)=\int_0^\omega L\big(x(t),\phi(t),\dot x(t),\dot \phi(t)\big)dt$$ attains a minimum in a set $E_{k,\omega}=X_\omega\times \Phi_{k,\omega}$. 

\begin{lemma}\label{dfgd}For any $(x,\phi)\in H^1_o(\mathbb{R})\times H^1_o(\mathbb{R})$ the following inequality holds
$$S(x,\phi)\ge \frac{1}{2}\|\dot\phi\|^2_{L^2(0,\omega)}+\frac{1}{2}\|\dot x\|^2_{L^2(0,\omega)} -\frac{\omega^{3/2}}{\sqrt 2}\|\dot x\|_{L^2(0,\omega)}. $$\end{lemma}
{\it Proof.} Indeed, with the help of Cauchy inequality it immediately follows that
\begin{align}S(x,\phi)&\ge\frac{1}{2}\|\dot\phi\|^2_{L^2(0,\omega)}+\frac{1}{2}\|\dot x\|^2_{L^2(0,\omega)}-\|x\|_{L^1(0,\omega)}\nonumber\\
&\ge \frac{1}{2}\|\dot\phi\|^2_{L^2(0,\omega)}+\frac{1}{2}\|\dot x\|^2_{L^2(0,\omega)}-\|x\|_{L^2(0,\omega)}\sqrt\omega.\nonumber
\end{align}
It remains to apply Lemma \ref{xddddd}.

The Lemma is proved.

\subsection{Minimization of the Action Functional}
Let $\{(x_n,\phi_n)\}_{n\in\mathbb{N}}\subset  E_{k,\omega}$ be a minimizing sequence for the functional $S$ that is 
$$S(x_n,\phi_n)\to \alpha,\quad n\to\infty,\quad \alpha=\inf_{E_{k,\omega}} S.$$
From Lemma \ref{dfgd} it follows that the sequence $\{(x_n,\phi_n)\}_{n\in\mathbb{N}}$ is bounded in $X_\omega\times H^1_o(\mathbb{R})$ and $\alpha>-\infty.$

Thus the sequence $\{(x_n,\phi_n)\}$ contains a weakly convergent subsequence, we denote this subsequence  by the same letters: $$x_n\to x_*\in X_\omega,\quad \phi_n\to \phi_*\in H_o^1(\mathbb{R}).$$ Since a  convex set of a locally convex space is closed iff it is weakly closed \cite{edv}, we have $\phi_*\in \Phi_{k,\omega}$.

We also know from analysis that the sequence $\{(x_n,\phi_n)\}$ contains a subsequence that is convergent in $C[0,\omega]\times C[0,\omega]$. (See Remark \ref{dfgdggg}.) So we accept that  $(x_n,\phi_n)\to (x_*,\phi_*)$ in $C[0,\omega]\times C[0,\omega]$.  

Our next goal is to prove that $\alpha=S(x_*,\phi_*).$

Observe the following evident estimates
\begin{align}
\int_0^\omega x_n^2\dot \phi_n^2dt&\ge \int_0^\omega(x_n^2-x_*^2)\dot \phi_n^2 dt
\nonumber\\&+ \int_0^\omega x_*^2\dot \phi_*^2dt +2\int_0^\omega\dot \phi_*x_*^2(\dot\phi_n-\dot \phi_*)dt,\label{qdgq}\\
\int_0^\omega \dot \phi_n^2dt&\ge 
 \int_0^\omega \dot \phi_*^2dt +2\int_0^\omega\dot \phi_*(\dot\phi_n-\dot \phi_*)dt.\label{qsdfsdfqq}\end{align}
Since $x_n\to x_*$ in $C[0,\omega]$ and the sequence $\{\dot \phi_n\}$ is bounded in $L^2(0,\omega)$ the first term in the right side of formula (\ref{qdgq}) vanishes as $n\to \infty.$

The last terms in the right sides of formulas (\ref{qdgq}) and (\ref{qsdfsdfqq}) are vanished as $n\to \infty$ because $\phi_n\to \phi_*$ weakly in $H_o^1(\mathbb{R})$.

Note also that $$\int_0^\omega x_n\sin\phi_ndt\to \int_0^\omega x_*\sin\phi_*dt$$ it is because $\{(x_n,\phi_n)\}$ converges in $C[0,\omega]\times C[0,\omega]$. 

Gathering all these observations we get $\alpha\ge S(x_*,\phi_*)$. So that  $$\alpha=S(x_*,\phi_*),\quad (x_*,\phi_*)\in E_{k,\omega}.$$ 

\subsection{Weak Solutions to the Lagrange Equations}
Take any two functions $x,\phi\in X_\omega$ and put
$$f(\xi,\eta)=S(x_*+\xi x,\phi_*+\eta \phi),\quad \xi,\eta\in \mathbb{R}.$$
From previous section it follows that a point $\xi=\eta=0$ is a minimum of $f$. This implies
$$\frac{\partial f}{\partial\xi}\Big|_{\xi=\eta=0}=\frac{\partial f}{\partial\eta}\Big|_{\xi=\eta=0}=0,$$
or in the detailed form
\begin{align}
\int_0^\omega \Big(\dot x_*(t)\dot x(t)dt&+\frac{\partial L}{\partial x}\big(x_*(t),\phi_*(t),\dot x_*(t),\dot \phi_*(t)\big)x(t)\Big)dt=0,\label{sdge5}\\
\int_0^\omega\Big( (1+x_*^2(t))\dot\phi_*(t)\dot\phi(t)dt&+\frac{\partial L}{\partial \phi}\big(x_*(t),\phi_*(t),\dot x_*(t),\dot \phi_*(t)\big)\phi(t)\Big)dt=0.\label{srt42}
\end{align}
Equations (\ref{sdge5}) and (\ref{srt42}) imply that the functions $x_*,\phi_*$ are the weak solutions to the Lagrange equations and  $x,\phi\in X_\omega$ are the test functions.

\subsection{Regularization}

From the theory developed above we know that $x_*,\phi_*$ belong to $H^1_{\mathrm{loc}}(\mathbb{R})$ end by the Sobolev embedding theorem $x_*,\phi_*\in C(\mathbb{R}).$

Our aim is to show that $x_*,\phi_*\in C^2(\mathbb{R}).$ Let us check this for $\phi_*$;  the corresponding result for $x_*$ follows in the same way.

Introduce a space
$$Y_\omega=\{u\in L^2_{\mathrm{loc}}(\mathbb{R})\mid u(-t)=u(t),\quad u(t+\omega)=u(t)\}.$$
In this definition the equalities hold almost everywhere.

Assume that a function $y$ belongs to $Y_\omega$. If in addition this function satisfies equality 
$$\int_0^\omega y(s)ds=0$$ then $$\phi(t)=\int_0^ty(s)ds\in X_\omega.$$ Moreover, it is clear that every function from $X_\omega$ can be presented in this way. 

Let us put
$$ a(t)=(1+x_*^2(t))\dot\phi_*(t)\in Y_\omega,\quad l(t)=\frac{\partial L}{\partial \phi}\big(x_*(t),\phi_*(t),\dot x_*(t),\dot \phi_*(t)\big)\in X_\omega.$$
Introduce the following linear functionals
\begin{align}p(y)&=\int_0^\omega y(s)ds,\nonumber\\ 
h(y)&=\int_0^\omega\Big(a(t)y(t)+l(t)\int_0^ty(s)ds\Big)dt.\nonumber
\end{align}
Them both belong to $Y'_\omega$.
By  Fubini's theorem we can rewrite the last functional in the form
$$h(y)=\int_0^\omega a(t)y(t)dt+\int_0^\omega y(s) \int_s^\omega l(t)  dtds.$$
From equation (\ref{srt42}) we know that $\ker p\subset\ker h.$ Therefore there exists a constant $\lambda$ such that \begin{equation}\label{ddfg55}h=\lambda p.\end{equation}
Since $y\in Y_\omega$ is an arbitrary function, and the functions $a(t),\,\int_t^\omega l(s)ds$ are even,
 equation (\ref{ddfg55}) takes the form
$$(1+x_*^2(t))\dot\phi_*(t)+\int_t^\omega l(s)ds=\lambda.$$
Since $\{l,x_*\}\subset X_\omega\subset C(\mathbb{R})$ we obtain $\phi_*\in C^2(\mathbb{R})$.

The Theorem is proved.



\begin{thebibliography}{99}

\bibitem{edv} R. Edwards: Functional Analysis. New York, 1965.
\bibitem{Shefer} H.  Schaefer:  Topological Vector Spaces, London, 1966.
\bibitem{iosida}  K. Yosida: Functional Analysis. Springer, Berlin, 1965.

\end{thebibliography}
\end{document}